\documentclass{amsart}
\usepackage[utf8]{inputenc}
\usepackage{bbold}
\usepackage{amsmath, amssymb, comment}
\usepackage{graphicx, hyperref}
\usepackage{wrapfig}
\usepackage{amsthm}
\usepackage{color}
\usepackage{float}
\usepackage{tikz} 
\usetikzlibrary{cd, arrows.meta}
\usepackage{quiver}
\graphicspath{{./fig/}}
\usepackage[capitalize, nameinlink, noabbrev, nosort]{cleveref}

\usepackage[utf8]{inputenc}
\usepackage{amsmath, color}
\usepackage{graphicx}
\usepackage{wrapfig}
\usepackage{amsthm}
\usepackage[margin=1in]{geometry}

\definecolor{lightviolet}{RGB}{197,180,22}
\definecolor{darkcyan}{RGB}{0, 170, 170}

\DeclareMathOperator{\wt}{wt}

\numberwithin{equation}{section}

\newtheorem{theorem}{Theorem}[section]

\newtheorem{lemma}[theorem]{Lemma}

\newtheorem{conjecture}[theorem]{Conjecture}

\theoremstyle{definition}
\newtheorem{definition}[theorem]{Definition}

\newtheorem{notation}[theorem]{Notation}
\newtheorem{convention}[theorem]{Convention}
\newtheorem{remark}[theorem]{Remark}
\newtheorem{question}[theorem]{Question}

\def\tDelta{\widetilde{\Delta}}
\def\tone{T1~}
\def\ttwo{T2~}
\def\G{\mathcal{G}}

\title{A dimer view on Fox's trapezoidal conjecture}
\author{Karola M\'esz\'aros, Melissa Sherman-Bennett, Alexander Vidinas}

\begin{document}
\maketitle
\begin{abstract}
     Fox's conjecture (1962) states that the sequence of absolute values of the coefficients of the Alexander polynomial of alternating links is trapezoidal. While the conjecture remains open in general, a number of special cases have been settled, some quite recently:  Fox's conjecture was shown to hold for special alternating links by Hafner, M\'esz\'aros, and Vidinas (2023) and for certain diagrammatic Murasugi sums of special alternating links by Azarpendar,  Juh\'asz, and K\'alm\'an (2024).  In this paper, we give an alternative proof of  Azarpendar, Juh\'asz, and K\'alm\'an's aforementioned beautiful result via a dimer model for the Alexander polynomial. In doing so, we not only obtain a significantly shorter proof of  Azarpendar,  Juh\'asz, and K\'alm\'an's result than the original, but we also obtain several theorems of independent interest regarding the Alexander polynomial, which are readily visible from the dimer point of view.
\end{abstract}

\section{Introduction}

The central question of knot theory is that of distinguishing links up to isotopy. The first polynomial invariant of links devised to help answer this question was the Alexander polynomial in 1928 \cite{alexander1928topological}. Fox famously conjectured in 1962 \cite{fox} that the sequence of absolute values of the coefficients of the Alexander polynomial of alternating links is trapezoidal. Denote the Alexander polynomial of an oriented link $L$ by $\Delta_L(t)$. Recall that a sequence $a_0, a_1, \ldots, a_n$ is trapezoidal if: $a_0<a_1<\cdots <a_k=\cdots=a_m>a_{m+1}>\cdots > a_n$  for some $0\leq k\leq m\leq n$.

\begin{conjecture}[\cite{fox}] 
        Let $L$ be an alternating link. Then the absolute values of the coefficients of the Alexander polynomial $\Delta_L(t)$ form a trapezoidal  sequence.
\end{conjecture}

Fox's conjecture remains open in general, although some special cases have been settled by Hartley \cite{H79} for two-bridge knots,  Murasugi \cite{murasugi} for a family of alternating algebraic links, and Ozsv\'ath and Szab\'o \cite{OS03} for the case of genus $2$ alternating knots, 
 among others.  That Fox's conjecture holds for genus $2$ alternating knots was also confirmed by Jong 
\cite{jong2009alexander}. 

More recently, Fox's conjecture was settled for special alternating links by Hafner, M\'esz\'aros and Vidinas  \cite{HMV-paper1}, and an alternative proof of the same result was given by K\'alm\'an, M\'esz\'aros and Postnikov \cite{kmp2025}. We recall that a \emph{special alternating link} admits a diagram which is alternating and has only type 1 Seifert circles\footnote{A Seifert circle is type 1 if it bounds a region of the link diagram, and otherwise is type 2.}. Using this result as a base case, with a beautiful insight, Azarpendar, Juh\'asz and K\'alm\'an \cite{kalman} settled Fox's conjecture for alternating links which admit an alternating diagram whose type 2 Seifert circles have small ``length" (see Definition \ref{def:length}). We note that \cite{kalman} use the language of ``diagrammatic Murasugi sum" to describe the class of alternating links in their result.

The main object of this paper is to (1) translate  Azarpendar, Juh\'asz and K\'alm\'an's aforementioned result into the dimer language, (2) give a simpler proof of it than the original, and (3) unravel why the main idea of the proof currently only works when the type 2 Seifert circles are small in length (in particular, length at most 2). If this restriction on length could be removed, then Fox's conjecture -- after quite some time -- would finally be settled.

To obtain our alternative proof, we use a dimer model for the Alexander polynomial similar to those used in \cite{cohen_dissertation,TwistedDimer,perfect_matchings,cluster_alg}. Utilizing this model, we give a short proof of the following beautiful theorem originally discovered and proved by Azarpendar, Juh\'asz and K\'alm\'an \cite{kalman}. See Sections \ref{sec:alexpoly}, \ref{sec:SA} and \ref{sec:setup} for the relevant definitions.

\begin{theorem}\cite[Theorem 2.23]{kalman} \label{main}
Let $L$ be an alternating link diagram whose type 2 Seifert circles have length at most 2. Then the absolute values of the coefficients of $\Delta_L(t)$ form a trapezoidal sequence.
\end{theorem}

We note that due to the technicalities involved, in \cite{kalman} only the statement and a sketch of the proof idea were given for the above theorem; a full proof was given for a special case. Our dimer model approach swiftly takes care of the above theorem in its stated generality, with many -- if not all! -- of the technicalities gone. It appears as Theorem \ref{thm:full_decomp} in the text below.

We also explain through computational evidence why the basic idea of the proof of Theorem \ref{main} does not seem to apply when some type 2 Seifert circle has length 3 or greater. We hope that identifying these apparent pitfalls will inspire others to find a way around them.

In addition to the above major goals, the dimer model also sheds further light on the Alexander polynomial, which we collect in a series of theorems. If $C$ is a type 2 Seifert circle of diagram $L$, one can modify the crossings along $C$ to obtain two new diagrams $L', L''$ (see Definition \ref{def:sum_on_circle} for details) with fewer type 2 Seifert circles. We use the notation $L=L' *_C L''$, and, following \cite{kalman}, say that $L$ is the \emph{diagrammatic Murasugi sum} of $L'$ and $L''$. 

The next theorem gives lower bounds on the absolute values of coefficients of $\Delta_L(t)$ when $L$ is alternating. For a polynomial $f$, we use the notation $|f|$ to mean ``replace all coefficients with their absolute values." Additionally, $\tDelta_L$ refers to the \textit{symmetrized Alexander polynomial}: an expression of the Alexander polynomial as a Laurent polynomial in variables $t^{1/2},t^{-1/2}$ whose degree sequence is centered around zero. 

\begin{theorem}
    Suppose $L$ is alternating and $L= L' *_C L''$.
    \begin{itemize}
        \item If $C$ is length $1$, then $|\tDelta_L| =  |\tDelta_{L'}||\tDelta_{L''}|$;
        \item If $C$ is length $l \geq 2$, then 
        $$|\tDelta_L| \geq  |\tDelta_{L'}||\tDelta_{L''}| + |\tDelta_{\widetilde{L'}}||\tDelta_{\widetilde{L''}}| $$ coefficient-wise, where $\widetilde{L'}, \widetilde{L''}$ are obtained from $L', L''$ as described in Definition \ref{def:y-moves-links}.
    \end{itemize}
\end{theorem}

The above theorem appears as Theorem \ref{thm:alternating-sums-non-canceling} in the text. The length $1$ and $2$ cases of the above theorem also appear in \cite{kalman}.

\begin{theorem} 
     Suppose $L$ is an alternating link diagram and $L = L' *_C L''$. Then $\tDelta_L$ and $\tDelta_{L'}\tDelta_{L''}$ have the same support.    
\end{theorem}
Above,  the \textbf{support} of  $f\in\mathbb{Z}[t^{1/2},t^{-1/2}]$ is the subset $S$ of   $\{z/2 \ | \ z\in\mathbb{Z}\}$ such that the monomial $t^s$ appears in $f$ with nonzero coefficient if and only if $s \in S$. The above theorem appears as Theorem \ref{thm:terms} in the text. 

\medskip

\noindent \textbf{Roadmap of the paper.} In Section \ref{sec:alexpoly} we lay out the combinatorial dimer model that we utilize in the paper. In Section \ref{sec:SA} we give an exposition of special alternating links, alternating links, and diagrammatic Murasugi sums. In Section \ref{sec:setup} we translate the diagrammatic Murasugi sum into dimer language. The latter leads us to the proofs of our main results, namely Theorems \ref{thm:term_comparison},\ref{thm:alternating-sums-non-canceling}, \ref{thm:terms} and  \ref{thm:full_decomp}  in Section \ref{sec:proof}. 
Finally, in Section \ref{sec:obstacle} we give computations showing why our approach did not extend  to all alternating links.

\section{The face-crossing incidence graph and Alexander polynomial of a link}
\label{sec:alexpoly}

In this section, we recall Kauffman's state summation formula for the (symmetrized) Alexander polynomial of a link. We then outline how to reinterpret states as perfect matchings on the \emph{truncated face-crossing incidence graph}, and the state sum as a sum over edge weights of dimers, as in prior work of Cohen \cite{TwistedDimer}.

\subsection{The symmetrized Alexander polynomial}\label{sec:symm_alex}

The Alexander polynomial of a link $L$ is a Laurent polynomial $\Delta_L \in \mathbb{Z}[t^{\pm 1/2}]$ which is defined up to multiplication by a signed Laurent monomial in $t^{1/2}$. More precisely, for $f, g \in \mathbb{Z}[t^{\pm 1/2}]$, define an equivalence relation $\sim$ by $f \sim g$ if $f= \pm t^{m/2} g$ for some integer $m$. Then $\Delta_L$ denotes a $\sim$-equivalence class, which is an isotopy invariant of the link $L$. The sequence of absolute values of the coefficients of $\Delta_L$ is palindromic (which follows from, e.g., \cite[pg. 183]{K06}\label{thm:skein_var}).

The \emph{symmetrized Alexander polynomial}\footnote{Also called the Conway normalization of the Alexander polynomial in \cite{K06}} $\tDelta_L \in \mathbb{Z}[t^{\pm 1/2}]$ is an element of the equivalence class $\Delta_L$ whose support is centered around zero. That is, if the highest degree term of $\tDelta_L$ is $t^a$, then the lowest degree term is $t^{-a}$. This determines $\tDelta_L$ up to sign; our conventions below will also fix the sign. Note that either all exponents of $\tDelta_L$ are in $\mathbb{Z}$ (if $\tDelta_L$ has odd degree) or all exponents of $\tDelta_L$ are in $\mathbb{Z}+ \frac{1}{2}$ (if $\tDelta_L$ has even degree). Moreover, if the link $L$ is alternating, the support of $\tDelta_L$ is saturated with respect to these constraints, meaning that if $b$ and $b+a$ are in the support, so are $b+1, \dots, b+a-1$. With this in mind, we make the following definition.

\begin{definition}
 A Laurent polynomial $f$ in $\mathbb{Z}[t^{\pm 1/2}]$ is \textbf{trapezoidal} if for some $b \in \frac{1}{2}\mathbb{Z}$,
 \[f= a_0 t^{b} + a_1 t^{b+1} + \dots + a_n t^{b+n} \] 
  and the absolute values $|a_0|,|a_1|, \dots, |a_n|$ form a trapezoidal sequence.
\end{definition} 

We now review a \emph{state summation} formula for the Alexander polynomial of an oriented link, given by Kauffman \cite{K06}. For the particular choice of state weights given here, this formula yields the symmetrized Alexander polynomial $\tDelta_L$.

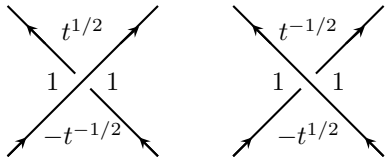
\begin{figure}[ht]
\centering
    \begin{tikzpicture}
        \draw[black, thick] (0,0)--(2,2); 
        \draw[black, thick] (0,2)--(0.9,1.1);
        \draw[black, thick] (2,0)--(1.1,0.9);
        \draw[-stealth, black, thick] (0,0) -- (0.25,0.25);
        \draw[-stealth, black, thick] (2,0) -- (1.75,0.25);
        \draw[-stealth, black, thick] (0.5,1.5) -- (0.25, 1.75);
        \draw[-stealth, black, thick] (1.5,1.5) -- (1.75,1.73);

        \draw[black, thick] (2+3,0)--(0+3,2); 
        \draw[black, thick] (0+3,0)--(0.9+3,0.9);
        \draw[black, thick] (1.1+3,1.1)--(2+3,2);
        \draw[-stealth, black, thick] (0+3,0) -- (0.25+3,0.25);
        \draw[-stealth, black, thick] (2+3,0) -- (1.75+3,0.25);
        \draw[-stealth, black, thick] (0.5+3,1.5) -- (0.25+3, 1.75);
        \draw[-stealth, black, thick] (1.5+3,1.5) -- (1.75+3,1.73);

        \node at (1,1.7) {$t^{1/2}$};
        \node at (1, 0.3) {$-t^{-1/2}$};
        \node at (1.4, 1) {$1$};
        \node at (0.6,1) {$1$};
        
        \node at (1+3,1.7) {$t^{-1/2}$};
        \node at (1+3, 0.3) {$-t^{1/2}$};
        \node at (1.4+3, 1) {$1$};
        \node at (0.6+3,1) {$1$};
    \end{tikzpicture}
\caption{The weighting of present regions around each crossing in the definition of the state sum.}
\label{fig:weighting}
\end{figure}

Fix an oriented link diagram\footnote{Abusing notation, we use the same symbol for an oriented link and a diagram for this link.} for $L$ and a distinguished \textbf{segment} $i$ of $L$ (that is, a curve between two crossings). The two planar regions adjacent to $i$ are called \textbf{absent}; the remaining regions are \textbf{present}. Weight the present regions around each crossing as in \ref{fig:weighting}.

A \textbf{state} is a bijection between crossings in the diagram and the present regions so that each crossing is mapped to one of the four regions it meets. This bijection is represented diagrammatically by adding a \textbf{state marker} to the region each corner is mapped to. We denote by $\mathcal{S}_{L,i}$ the set of states of the link diagram $L$ with distinguished segment $i$. The \textbf{weight} of a state $S$, denoted $\langle L|S\rangle$, is the product of the weights of the marked regions of $S$.

We take the following result of Kauffman   $\tDelta_L$ as our definition of $\tDelta_L$:

\begin{theorem}[{\cite[pg. 182]{K06}}]\label{thm:state-sum-alexander-poly}
    Let $L$ be an oriented link diagram and $i$ a segment of $L$. Then the symmetrized Alexander polynomial of $L$ is given by \[\tDelta_L=\sum_{S \in \mathcal{S}_{L,i}}\langle L|S\rangle.\]
\end{theorem}

\begin{remark}
    The symmetrized Alexander polynomial $\tDelta_L$ is independent of the choice of segment used to compute it. In \cite{K06}, the weights have no negative signs, and instead, signs are accounted for by the notion of a {\em black hole}; our signs exactly record the presence of a black hole. In \cite{kalman}, the authors' weighting differs from that of Figure \ref{fig:weighting} by exchanging $t^{1/2}\leftrightarrow -t^{1/2}$ and $t^{-1/2} \leftrightarrow -t^{-1/2}$. As a result, the polynomial computed by their method, when that polynomial has an even number of terms, differs from ours by multiplying the entire polynomial by $-1$. When that polynomial has an odd number of terms, it is identical to ours. 
\end{remark}

\subsection{The face-crossing incidence graph of a link diagram}

In this subsection, we define the \textit{truncated face-crossing incidence graph}\footnote{Also called the balanced overlaid Tait graph in \cite{TwistedDimer,cohen_dissertation,perfect_matchings}} $\G_{L,i}$ of a link diagram. We describe how to assign edge weights to $\G_{L,i}$ so the state sum of Theorem \ref{thm:state-sum-alexander-poly} becomes a weighted sum of perfect matchings on $G_{L,i}$. This reformulation of weighted states as weighted perfect matchings is due to Cohen \cite{TwistedDimer}, though we use the edge weights in Figure \ref{fig:weighting} rather than those used in \cite{TwistedDimer}.

\begin{definition}
    Let $L$ be a link diagram. Stereographically project $S^2$ through any point in the complement of the diagram. The \textbf{face-crossing incidence graph}, $\G_L$, is a plane graph defined as follows. Place a black vertex $b_c$ on each crossing $c$ of $L$ and a white vertex $w_r$ in each of its planar regions $r$. There is an edge $(b_c, w_r)$ in $G_L$ if and only if crossing $c$ touches region $r$. We sometimes abuse notation and identify $b_c$ with $c$ and $w_r$ with $r$.
\end{definition}

\begin{figure}
    \centering
    \includegraphics[width=0.7\linewidth]{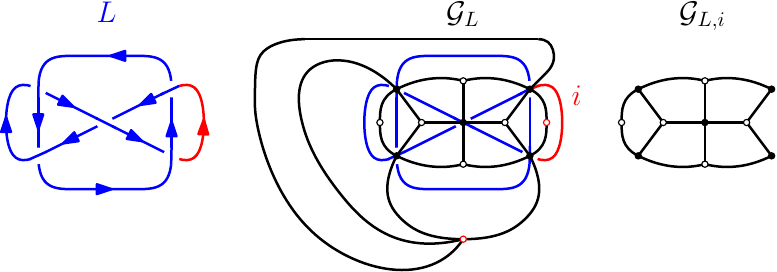}
    \caption{Left: In blue, a link diagram for the Whitehead link, with a choice of distinguished segment $i$ in red. Center: The face-crossing incidence graph $\G_L$. Right: The truncated face-crossing incidence graph $\G_{L,i}$.}
    \label{fig:whitehead}
\end{figure}

See Figure \ref{fig:whitehead} for an example of a face-crossing incidence graph.

For any fixed segment $i$ of $L$, stereographically projecting through any point in a region adjacent to $i$ yields a drawing of $L$ in the plane for which $i$ bounds the exterior region.

\begin{definition} \label{def:truncated}
    Let $L$ be a link diagram, and fix a segment $i$ of $L$. Changing the stereographic projection of $L$ if necessary, we may assume that $i$ bounds the exterior region. The \textbf{truncated face-crossing incidence graph} is the plane graph $\G_{L,i}$ obtained by deleting from $G_L$ the two white vertices of the infinite face along with any incident edges. 
\end{definition}

See Figure \ref{fig:whitehead} for examples of truncated face-crossing incidence graphs. To compute the Alexander polynomial of a link diagram from a truncated face-crossing incidence graph, we weight its edges according to the following convention.

\begin{definition}\label{def:edge_weight}
    Let $L$ be a link diagram and let $i$ denote a distinguished segment of $L$. Each edge $e$ of $\G_{L,i}$ receives a \textbf{weight} $\wt(e)$ as below.  
    
    \begin{center}
    \begin{tikzpicture}
        \draw[black, thick] (1,0) -- (1,.9);
        \draw[-stealth, black, thick] (1,1.1) -- (1,2);
        \draw[-stealth, black, thick] (0,1) -- (2,1);
        \node at (0.75+4,1.6) {$1$};
        \node at (0+4,.7) {$-t^{-1/2}$};
        \node at (2.05+4,1.5) {$t^{1/2}$};
        \node at (1.3+4,.4) {$1$};

        \draw[black, thick] (4,0) -- (5,1);
        \draw[black, thick] (4,2) -- (5,1);
        \draw[black, thick] (6,0) -- (5,1);
        \draw[black, thick] (6,2) -- (5,1);
        \fill (5,1) circle (2pt);
        \filldraw[color=black, fill=white, very thick] (6,2) circle (2pt);
        \filldraw[color=black, fill=white, very thick] (6,0) circle (2pt);
        \filldraw[color=black, fill=white, very thick] (4,2) circle (2pt);
        \filldraw[color=black, fill=white, very thick] (4,0) circle (2pt);
\end{tikzpicture}
\hspace{1in}
\begin{tikzpicture}

        \draw[-stealth, black, thick] (1,-3) -- (1,-1);
        \draw[black, thick] (0,-2) -- (0.9,-2);
        \draw[-stealth, black, thick] (1.1,-2) -- (2,-2);
        \node at (0.75+4,1.6-3) {$1$};
        \node at (0+4,.7-3) {$-t^{1/2}$};
        \node at (2.05+4,1.5-3) {$t^{-1/2}$};
        \node at (1.3+4,.4-3) {$1$};

        \draw[black, thick] (4,0-3) -- (5,1-3);
        \draw[black, thick] (4,2-3) -- (5,1-3);
        \draw[black, thick] (6,0-3) -- (5,1-3);
        \draw[black, thick] (6,2-3) -- (5,1-3);
        \fill (5,1-3) circle (2pt);
        \filldraw[color=black, fill=white, very thick] (6,2-3) circle (2pt);
        \filldraw[color=black, fill=white, very thick] (6,0-3) circle (2pt);
        \filldraw[color=black, fill=white, very thick] (4,2-3) circle (2pt);
        \filldraw[color=black, fill=white, very thick] (4,0-3) circle (2pt);
    \end{tikzpicture}
    \end{center}
\end{definition}

\begin{figure}
    \centering
    \includegraphics[width=0.4\linewidth]{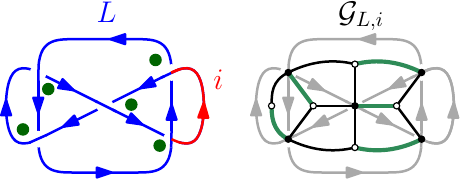}
    \caption{A diagram of the Whitehead link. Left: A Kauffman state on the link diagram. The distinguished segment $i$ is in red. Right: The corresponding matching on $\G_{L,i}$, with the same distinguished segment chosen.}
    \label{fig:whitehead_comparison}
\end{figure}

    Recall that, for a link diagram $L$ and choice of a segment $i$, a \emph{Kauffman state} is a bijection between present regions and crossings of a link diagram. Also note the vertices of $\G_{L,i}$ are precisely the crossings and present regions of $L$, and the edges connect each region to the crossings it meets. The next lemma gives a bijection between states and perfect matchings of $\G_{L,i}$; see Figure \ref{fig:whitehead_comparison} for an example. This lemma in a slightly different form can be found in \cite{perfect_matchings, cluster_alg}.

\begin{lemma}\label{rem:state_and_matching}
Let $L$ be a link diagram and let $i$ be a segment of $L$. Fix an embedding of $L$ for which $i$ is on the boundary of the exterior region. Let $S\in\mathcal{S}_{L,i}$, and define the perfect matching $M_S$ on ${\G_{L,i}}$ by $e=(b_c,w_r)\in M_S$ if the state marker near crossing $c$ in $S$ is in region $r$ (see below). The map $S \mapsto M_S$ is a weight-preserving bijection between states on $L$ with chosen segment $i$ and perfect matchings on $\G_{L,i}$.
\end{lemma}

    \begin{center}
    \begin{tikzpicture}
        \draw[black, thick] (1,0) -- (1,2);
      
        \draw[black, thick] (0,1) -- (2,1);
        \filldraw[color=blue] (1.2,1.2) circle (2pt);

        \draw[black, thick] (4,0) -- (5,1);
        \draw[black, thick] (4,2) -- (5,1);
        \draw[black, thick] (6,0) -- (5,1);
        \draw[blue, very thick] (6,2) -- (5,1);
        \fill (5,1) circle (2pt);
        \filldraw[color=black, fill=white, very thick] (6,2) circle (2pt);
        \filldraw[color=black, fill=white, very thick] (6,0) circle (2pt);
        \filldraw[color=black, fill=white, very thick] (4,2) circle (2pt);
        \filldraw[color=black, fill=white, very thick] (4,0) circle (2pt);
    \end{tikzpicture}
    \end{center}

\begin{notation}
In light of Remark \ref{rem:state_and_matching}, we abuse notation and, for the remainder of the paper, use $\mathcal{S}_{L,i}$ to refer to the set of perfect matchings of $\G_{L,i}$.
\end{notation}

As an immediate corollary of the above lemma, we can compute $\tDelta_L$ using matchings of $G_{L,i}$.

\begin{theorem}\label{thm:dimer_version}
    Let $L$ be a link diagram. The symmetrized Alexander polynomial can be computed as

    \[\tDelta_L = \sum_{M\in\mathcal{S}_{L,i}}\prod_{e\in M}\wt(e) =\sum_{M\in\mathcal{S}_{L,i}}\wt(M) .\]
\end{theorem}

\section{Special alternating links, alternating links, and the diagrammatic Murasugi sum}
\label{sec:SA}

In this section, we define  
\textit{special alternating links} and the \textit{diagrammatic Murasugi sum} of special alternating links following  \cite[Remark 2.17]{kalman}. A key property of  any alternating link is that it is a diagrammatic Murasugi sum of  special alternating links \cite{murasugi1965}.

\subsection{Seifert circles and special alternating links}
\begin{figure}[ht]
    \centering
    \includegraphics[width=0.5\linewidth]{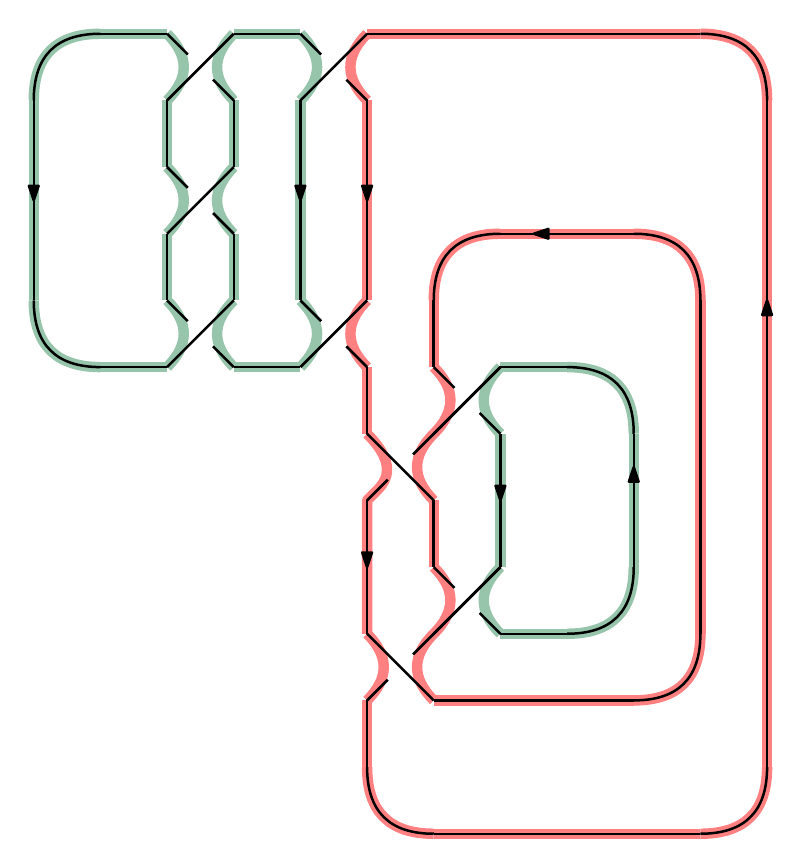}
    \caption{A link diagram and its Seifert circles. The \tone circles are shaded green and the \ttwo   circles are shaded red. }
    \label{fig:seif_circles} 
\end{figure}

Let $L$ be a link diagram. By smoothing each crossing as below\footnote{The smoothing does not depend on over vs. under crossing, so we omit this information in the figure.}, one obtains a collection of (possibly nested) oriented circles.
    \begin{center}
    \begin{tikzpicture}
        \draw[black, thick] (0,0)--(2,2); 
        \draw[black, thick] (0,2)--(2,0);
        \draw[-stealth, black, thick] (0,0) -- (0.25,0.25);
        \draw[-stealth, black, thick] (2,0) -- (1.75,0.25);
        \draw[-stealth, black, thick] (0.5,1.5) -- (0.25, 1.75);
        \draw[-stealth, black, thick] (1.5,1.5) -- (1.75,1.73);

        \node at (3,1) {$\mapsto$};

        \draw[black, thick] (0+4,0) -- (0.5+4,0.5);
        \draw[black, thick] (2+4,0) -- (1.5+4,0.5);
        \draw[black, thick] (0+4,2) -- (0.5+4,1.5);
        \draw[black, thick] (2+4,2) -- (1.5+4,1.5);

        \draw[-stealth, black, thick] (0+4,0) -- (0.25+4,0.25);
        \draw[-stealth, black, thick] (2+4,0) -- (1.75+4,0.25);
        \draw[-stealth, black, thick] (0.5+4,1.5) -- (0.25+4, 1.75);
        \draw[-stealth, black, thick] (1.5+4,1.5) -- (1.75+4,1.73);

        \draw[black, thick] plot [smooth] coordinates {(0.5+4,0.5) (0.8+4,0.8) (0.8+4, 1.2) (0.5+4,1.5)};

        \draw[black, thick] plot [smooth] coordinates {(1.5+4,0.5) (1.2+4,0.8) (1.2+4, 1.2) (1.5+4,1.5)};
    \end{tikzpicture}
    \end{center}
 Each oriented circle specifies a collection of segments in the \textit{original} link diagram $L$. Each such collection of segments is called a \textbf{Seifert circle}. See Figure \ref{fig:seif_circles} for an example.

Recall that a \emph{region} of a link diagram $L$ is a connected component of the complement $S^2 \setminus L$. As in \cite{kalman} we will refer to type 1 and type 2 Seifert circles.

\begin{definition} 
    A Seifert circle of link diagram $L$ is \textbf{type $1$} if it bounds a region of $L$ and \textbf{type $2$} otherwise. We use the abbreviations \tone circle and \ttwo circle for ``type 1 Seifert circle" and ``type 2 Seifert circle."
\end{definition}

A link diagram is \textbf{special}\footnote{This terminology comes from Crowell \cite{crowell}, where type 1 Siefert circles are called \emph{special Seifert circuits}.} if all of its Seifert circles are type $1$. Every link admits a special diagram \cite[Proposition 13.14]{burde2002knots}. A link diagram is \textbf{alternating} if, tracing along any strand, it alternates between going over and under at crossings. A link is \textbf{special alternating} if it admits a diagram that is \emph{both} alternating and special, which we call a special alternating diagram.

\subsection{The diagrammatic Murasugi sum}

In this section, we define the \textit{diagrammatic Murasugi sum} of two link diagrams \textit{along a \ttwo circle}, following \cite[Remark 2.17]{kalman}. A key result of Murasugi is that   any alternating link is a diagrammatic Murasugi sum of  special alternating links  Murasugi \cite{murasugi1965}; we use the terminology diagrammatic Murasugi sum  as in \cite{kalman}. 

\begin{definition}\label{def:sum_on_circle}
    Let $L$ be an oriented link diagram that is not special. Let $C$ be a \ttwo circle of $L$. Observe that the complement $S^2 \setminus C$ of $C$ in the sphere has two connected components, $R'$ and $R''$. We say that a crossing $c$ of $L$ \textbf{belongs to} $R'$ (resp. $R''$) if at least two of the regions adjacent to it  are contained in $R'$ (resp. $R''$). 
    Let $L'$ denote the link diagram obtained by smoothing the crossings of $C$ that do not belong to $R'$, then removing the portion of the resulting diagram which is contained in $R''$ (a union of unlinked unknotted circles). Define $L''$ analogously. 
We say that $L$ is the \textbf{diagrammatic Murasugi sum of $L'$ and $L''$ along $C$}, denoted by $L = L'*_C L''$.
\end{definition}

\begin{figure}[ht]
    \centering
    \includegraphics[width=\linewidth]{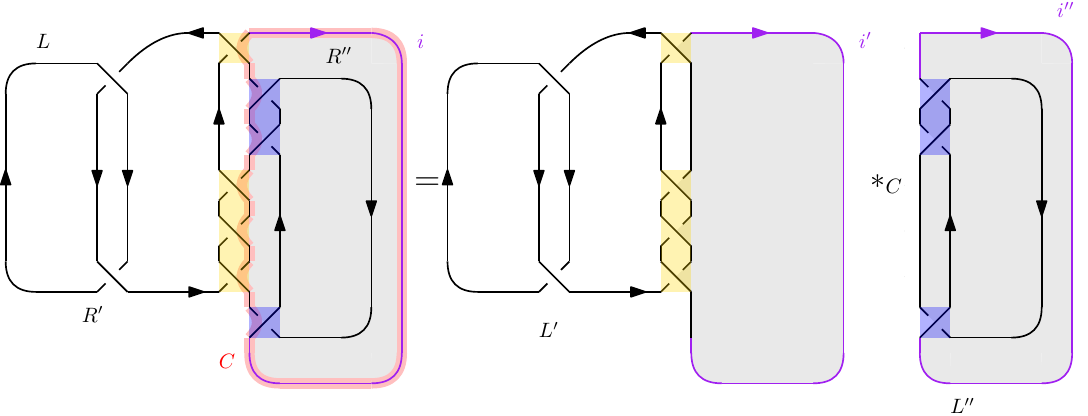}
    \caption{An illustration of diagrammatic Murasugi sum using the notation of Definition \ref{def:sum_on_circle}. On the left, a link diagram $L$ with a \ttwo circle $C$ highlighted in red and $R''$ shaded in grey. The crossings on $C$ belonging to $R'$ are in yellow, and those belonging to $R''$ are in blue. On the right, the links $L'$ and $L''$ so that $L = L' *_C L''$. The segment $i$ illustrates Convention \ref{conv:projection}, and $i', i''$ are the unique segments of $L', L''$ containing $i$, and are indicated in violet.}
    \label{fig:diag_sum_example}
\end{figure}

See Figure \ref{fig:diag_sum_example} for an example of Definition \ref{def:sum_on_circle}.

\begin{remark}\label{rem:decomp_tone}
We note that $C$ is a \tone circle for $L'$ and $L''$. Also, each \ttwo circle $C' \neq C$ of $L$ is a \ttwo circle of exactly one of $L'$ and $L''$, and all \ttwo circles of $L', L''$ arise in this way. As $L', L''$  both have strictly fewer \ttwo circles than $L$, we may iterate Definition \ref{def:sum_on_circle} and eventually obtain a collection of special diagrams. We also note that if $L$ is alternating, then $L', L''$ are alternating (cf. Lemma \ref{lem:y-moves-alternating-and-type-2-circles}), so iterating Definition \ref{def:sum_on_circle} on an alternating diagram produces a collection of special alternating diagrams.
\end{remark}

\section{Notation and Setup}
\label{sec:setup}
 
In this section we study the structure of truncated face-crossing incidence graphs of link diagrams which are diagrammatic Murasugi sums $L= L' *_C L''$. We    define the length of a diagrammatic sum, following the exposition of \cite{kalman}, and also introduce the closely related length of a \ttwo circle. Next we establish useful conventions in drawing the diagrammatic Murasugi sums we work with, as well as define the auxiliary diagrams $\widetilde{L'}$ and $\widetilde{L''}$ needed to prove the main results of the paper (Theorems \ref{thm:term_comparison},\ref{thm:alternating-sums-non-canceling}, \ref{thm:terms} and  \ref{thm:full_decomp}).

We begin with the defintion of the length of a \ttwo circle and of a diagrammatic Murasugi sum.

\begin{definition}\label{def:length}
Let $C$ be a \ttwo circle of $L$. Each segment of $C$ goes between two crossings, and each crossing belongs to either $R'$ or $R''$. The \textbf{length} of $C$ is defined to be
\[l=\frac{1}{2}|\{\text{segments } i \text{ of }C~\text{which connect a crossing belonging to}~R'~\text{and a crossing belonging to}~R''\}|.\]
If $L=L' *_C L''$, we also call $l$ the \textbf{length} of the diagrammatic sum $L' *_C L''$.
\end{definition}

Recall from Definition \ref{def:truncated} the truncated face-crossing incidence graph $\G_{L,i}$.

\begin{figure}
    \centering
    \includegraphics[width=0.9\linewidth]{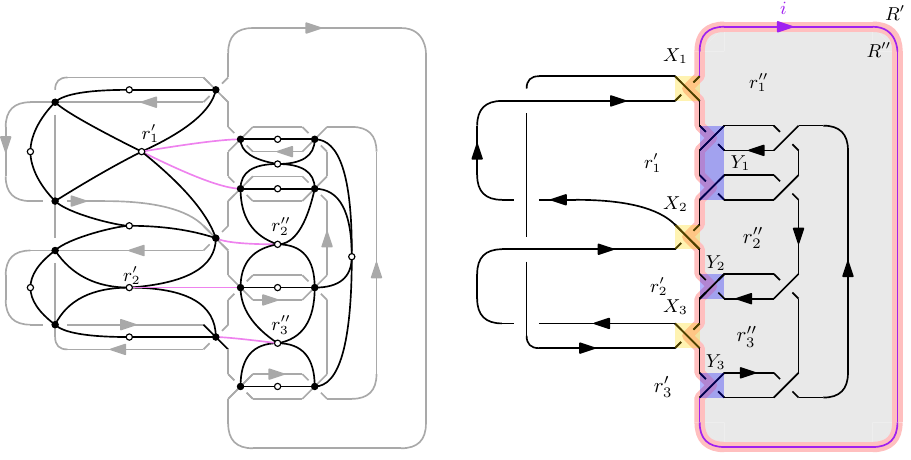}
    \caption{A summary of Notation \ref{not:sums}. Right: A link diagram. The \ttwo circle $C$ is indicated in red, and the region $R''$ bounded by $C$ is shaded grey. The groups of crossings $X_j$ are yellow, and the crossings $Y_j$ are blue. The regions $r'_j$ and $r''_j$ are labeled. Observe $r''_1$ and $r'_3$ are the absent regions from selecting the purple segment $i$.   Left: The graph $\G_{L,i}$ for the link diagram on the right, with its flock edges colored pink.}
    \label{fig:diagrammatic_sum_cartoon}
\end{figure}

\begin{convention}\label{conv:projection}
In what follows, we will always fix a link diagram $L \in S^2$ and a \ttwo circle $C$. We will also fix a stereographic projection of $L$ as follows: choose any segment $i$ of $C$ connecting a crossing belonging to $R'$ to a crossing belonging to $R''$ (or, if $C$ has length 0, choose any segment of $C$), and project through a point in $R'$ in a region adjacent to $i$. See Figure \ref{fig:diagrammatic_sum_cartoon}. We will also choose the same segment $i$ as the distinguished segment of $L$, and use this segment and the aforementioned projection to compute $\G_{L,i}$. We consider the graph $\G_{L,i}$ as an edge-weighted graph, with weights given by Definition \ref{def:edge_weight}. 
\end{convention}

Recall from Section \ref{sec:symm_alex} the notion of an \emph{absent region} of a diagram.

\begin{notation}\label{not:sums}  
    Suppose $L = L' *_C L''$, where $C$ has length $l>0$. Choose a segment $i$ of $C$ and a projection of $L$ as in Convention \ref{conv:projection}. See Figure \ref{fig:diagrammatic_sum_cartoon}.
    \begin{itemize}
        \item  Starting from a point $p$ on $i$ and tracing $C$, we see a collection $X_j$ of crossings belonging to $R'$ then a collection $Y_j$ of crossings belonging to $R''$, for $j=1, \dots, l$. That is, first we see the crossings in $X_1$, then those in $Y_1$, then those in $X_2$, etc. and after we see the crossings in $Y_l$, we return to $p$.
        \item The crossings in $X_j$ are all adjacent to a region $r''_j$ in $R''$, and the crossings in $Y_j$ are all adjacent to a region $r'_j$ in $R'$. Note that if $i$ is the distinguished segment, then $r''_1$ and $r'_l$ are the absent regions.
        \item In $\G_{L,i}$, there are edges from $r''_j$ to each crossing in $X_j$ and from $r'_j$ to each crossing in $Y_j$. 
        We refer to all such edges as \textbf{flock edges} of $\G_{L,i}$.
    \end{itemize}
    If $C$ has length 0, then all crossings on $C$ are adjacent to an absent region, so $\G_{L,i}$ is disconnected and has no flock edges.
\end{notation}

\begin{lemma}\label{lem:flock-edges-1} Suppose $L = L' *_C L''$. Choose $i$ as in Notation \ref{not:sums}. The flock edges of $\G_{L,i}$ have weight 1.
\end{lemma}

\begin{proof}
    The flock edges go from a crossing $c$ on the \ttwo circle $C$, to a region $r$ which is adjacent to $c$ and either immediately to the right of $C$ or immediately to the left. From the definition of Seifert circle and the weighting of Definition \ref{def:edge_weight}, we conclude these edges have weight 1.
\end{proof}

The next lemma shows that deleting the flock edges of $\G_{L,i}$ recovers the truncated face-crossing incidence graphs of $L'$ and $L''$.

\begin{lemma}\label{lem:d_move}
    Suppose $L = L' *_C L''$ and we have chosen segment $i$ in $C$ as in Convention \ref{conv:projection}. Let $i'$ and $i''$ be the unique segments of $L'$ and $L''$ which contain $i$. Delete the flock edges of $\G_{L,i}$. The resulting weighted graph is $\G_{L',i'} \sqcup \G_{L'',i''}$.
\end{lemma}

\begin{figure}
        \centering
        \includegraphics[width=\linewidth]{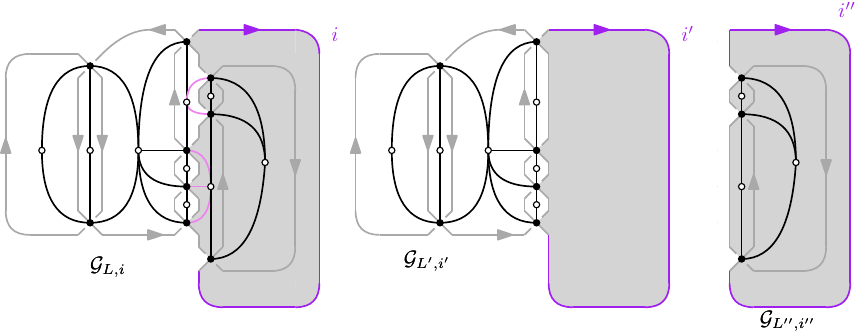}
        \caption{The graphs $\G_{L,i}$, $\G_{L', i'}$, and $\G_{L'', i''}$ for the sum $L= L' *_C L''$ in Figure \ref{fig:diag_sum_example}. The links are in grey, $R''$ is shaded, and the segments $i, i', i''$ are in violet. The flock edges of $\G_{L,i}$ are in pink. Deleting the flock edges produces $\G_{L', i'} \sqcup \G_{L'', i''}$.}

        \label{fig:d_move}
    \end{figure}

Refer to Figure \ref{fig:d_move} for an example.
    
\begin{proof}
We compare $\G_{L,i}$ with $\G_{L', i'}$; the comparison with $\G_{L'', i''}$ is similar. Recall that the crossings of $L'$ are exactly the crossings of $L$ which belong to $R'$. The non-infinite regions of $L'$ are exactly the regions of $L$ which are contained in $R'$, but the crossings in $\cup_j Y_j$ have been smoothed. So for non-infinite regions in $R''$, crossing incidence in $L'$ is exactly the same as crossing incidence in $L$ if you ignore the crossings $\cup_j Y_j$. Note also that the absent regions of $\G_{L,i}$ are $r_1''$ and $r_l'$; by inspection, the absent regions of $\G_{L',i'}$ are $r_l'$ and $R''$, which contains $r_1''$. 

Say a vertex of $\G_{L,i}$ is \emph{from} $R'$ if it is $b_c$ for $c$ a crossing belonging to $R'$ or $w_r$ for a region $r$ contained in $R'$; a vertex is from $R''$ in analogous circumstances. From the above observations, we conclude $\G_{L', i'}$ is the induced subgraph of $\G_{L,i}$ obtained by deleting all vertices from $R''$.

A similar argument shows that $\G_{L'', i''}$ is the induced subgraph of $\G_{L,i}$ obtained by deleting all vertices from $R'$.

 The flock edges are exactly the edges of $\G_{L,i}$ which connect a vertex from $R'$ to a vertex from $R''$. So $\G_{L', i'}$ and $\G_{L'', i''}$ are exactly the connected components of $\G_{L,i}$ with the flock edges removed.

To see the desired equality at the level of weighted graphs, we note that the weighting is determined by a local rule at each crossing. The neighborhood of every crossing of $L'$ and $L''$ is identical to that of some corresponding crossing of $L$, so the weights on nearby edges are also identical.

\end{proof}

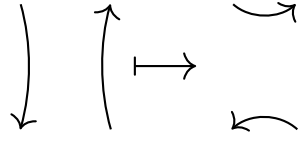
\begin{figure}
    \centering
\[\begin{tikzcd}
	{\textcolor{white}{\bullet}} & {\textcolor{white}{\bullet}} & {\textcolor{white}{\bullet}} & {\textcolor{white}{\bullet}} \\
	& {\textcolor{white}{\bullet}} & {\textcolor{white}{\bullet}} \\
	{\textcolor{white}{\bullet}} & {\textcolor{white}{\bullet}} & {\textcolor{white}{\bullet}} & {\textcolor{white}{\bullet}}
	\arrow[curve={height=-6pt}, from=1-1, to=3-1, thick]
	\arrow[curve={height=12pt}, from=1-3, to=1-4, thick]
	\arrow[maps to, from=2-2, to=2-3, thick]
	\arrow[curve={height=-6pt}, from=3-2, to=1-2, thick]
	\arrow[curve={height=-12pt}, tail reversed,no head, from=3-3, to=3-4, thick]
\end{tikzcd}\]
\caption{A swap move on a link diagram.}
\label{fig:swap-move}
\end{figure}

Next we will show that deleting certain non-flock edges of $\G_{L,i}$ produces the truncated face-crossing incidence graphs of two links. We will first define these links, and then state and prove the lemma. The definition of these links in the case where $C$ has length at most $2$ originally appears in \cite[proof of Theorem 2.20, Remark 2.24]{kalman}. Our definition will use the \textbf{swap move} demonstrated in Figure \ref{fig:swap-move}.

\begin{definition}\label{def:y-moves-links}
    Suppose $L = L' *_C L''$ has length at least 2. We use the notation from Notation \ref{not:sums}, and use $i', i''$ to denote the unique segments of $L', L''$ containing $i$. The link diagram $\widetilde{L'}$ is obtained from $L'$ by performing \textbf{swap moves} on $i'$ and the segment separating $X_j$ from $X_{j+1}$, for $j \in [l-1]$. One of the swap moves creates a segment between the first and last crossing in $X_1$, which we denote by $\widetilde{i'}$. The link diagram $\widetilde{L''}$ is obtained analogously from $L''$: first change the projection of $L''$ so that the other region adjacent to $i''$ is infinite, and then do swap moves on $i''$ and segments separating $Y_j$ and $Y_{j+1}$, for $j \in [l-1]$. One of these swap moves creates a segment between the first and last crossing in $Y_l$, which we denote by $\widetilde{i''}$.
\end{definition}

\begin{figure}
        \centering
        \includegraphics[width=0.6\linewidth]{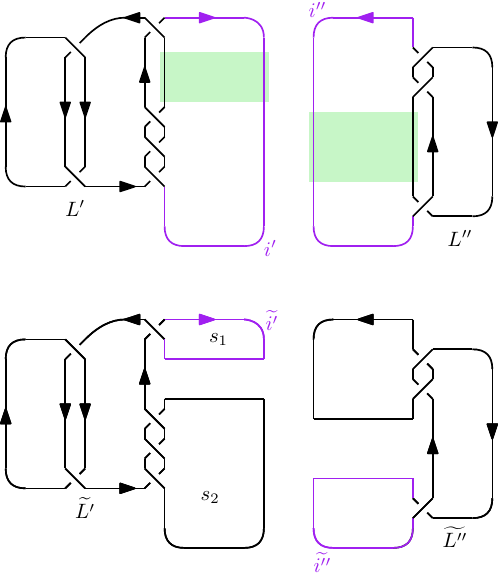}
        \caption{The construction of $\widetilde{L'}$ and $\widetilde{L''}$ from Definition \ref{def:y-moves-links}, using the Murasugi sum from Figure \ref{fig:diag_sum_example}, which has length 2. The area where the swap move is applied is shaded. The regions $s_1$ and $s_2$ referenced in the proof of Lemma \ref{lem:y_moves}.}
        \label{fig:y_move}
    \end{figure}

See Figure \ref{fig:y_move} for an example of $\widetilde{L'}$ and $\widetilde{L''}$.

\begin{lemma}\label{lem:y_moves}
    Suppose that $L=L' *_C L''$ has length at least 2. Let $\widetilde{L'}, \widetilde{L''}, \widetilde{i'}, \widetilde{i''}$ be as in Definition \ref{def:y-moves-links}. In $\G_{L,i}$, delete all non-flock edges adjacent to $r'_j, r''_j$ for $j \in [l].$ The resulting weighted graph is 
    \[\G_{\widetilde{L'},\widetilde{i'}} \sqcup \G_{\widetilde{L''},\widetilde{i''}}.\]
\end{lemma}

\begin{figure}
        \centering
        \includegraphics[width=0.7\linewidth]{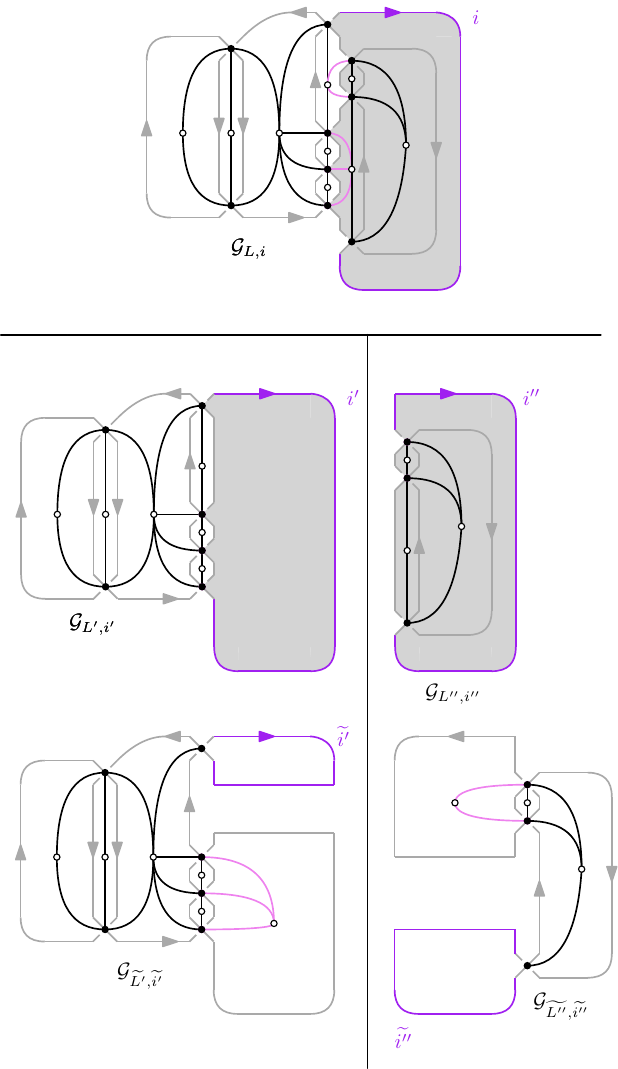}
        
        \caption{The graphs $\G_{L,i}$,  $\G_{{L'}, {i'}}$, $\G_{{L''}, {i''}}$, $\G_{\widetilde{L'}, \widetilde{i'}}$, and $\G_{\widetilde{L''}, \widetilde{i''}}$.
        }
        \label{fig:y_move_graphs}
    \end{figure}

\begin{proof} 
We consider how the diagram $\widetilde{L'}$ relates to the diagram $L'$. 
The swap moves applied to $L'$ replace the unique bounded region enclosed by $i'$ and create new regions $s_1, \dots, s_l$. See Figure \ref{fig:y_move} for an example. The region $s_j$ is adjacent to the crossings in $X_j$ and no other crossings for all $j \in [l]$. The swap move also identifies each region $r'_k$, where $k\in[l-1]$, with the infinite region. We choose the absent regions to be  $s_1$ and the infinite region. Thus, $\G_{\widetilde{L'}, \widetilde{i'}}$ is obtained from $\G_{L',i'}$ by adding new vertices ${s_j}$ where $j \in \{2, \ldots, l\}$, and edges $({s_j},c)$ for $c \in X_j$, $j \in \{2, \ldots, l\}$, then deleting the vertices $r'_k$, $k\in[l-1]$. By inspection, the new edges have weight $1$. Thus, $\G_{\widetilde{L'}, \widetilde{i'}}$ is isomorphic to the subgraph of $\G_{L,i}$ obtained by adding the flock edges $({r''_j},c)$ for $c \in X_j$, $j\in \{2,\dots l\}$, to $\G_{L',i'}$ and removing all edges incident to vertices $r'_k$, $k\in[l-1]$. See Figure \ref{fig:y_move_graphs} for an example. By Lemma \ref{lem:flock-edges-1}, the edges $({r''_j},c)$, $j\in\{2,\ldots,l\}$, also have weight $1$, so moreover, this isomorphism preserves edge weights.  

Similarly, $\G_{\widetilde{L''}, \widetilde{i''}}$ is obtained from $\G_{L'',i''}$ by adding new vertices ${t_j}$ and weight $1$ edges $({t_j}, c)$ for $c \in Y_j$, $j \in [l-1]$. Thus,  $\G_{\widetilde{L''}, \widetilde{i''}}$ is isomorphic to the subgraph of $\G_{L,i}$ obtained by adding the flock edges $({r'_j},c)$ for $c \in Y_j$, $j\in [ l-1]$, to $\G_{L'',i''}$ and removing all edges incident to vertices $r''_k$, $k\in\{2,\dots,l\}$. Moreover, this isomorphism preserves edge weights.  
\end{proof}

We will also need the following lemma. Related statements have also been developed in \cite{kalman}.

\begin{lemma}\label{lem:y-moves-alternating-and-type-2-circles}
    Suppose that $L=L' *_C L''$ has length at least 2 and let $\widetilde{L'}, \widetilde{L''}$ be as in Definition \ref{def:y-moves-links}. Then $L'$ and $\widetilde{L'}$ (resp. $L''$ and $\widetilde{L''}$) have the same number of \ttwo circles. If $L$ is alternating, then so are $L'$, $L''$, $\widetilde{L'}$, and $\widetilde{L''}$.
\end{lemma}

\begin{proof}
    For the first statement, we recall by Remark \ref{rem:decomp_tone} that $C$ is a \tone circle in each of $L'$ and $L''$. Applying the swap moves on $L'$ (or $L''$) will change the \tone circle $C$ into $l$ \tone circles; these are the circles enclosing the regions $s_1$ and $s_2$ on the left in Figure \ref{fig:y_move}. All other Seifert circles of $L'$ (or $L''$) are undisturbed by the swap moves.

    For the second statement, we focus on $L'$, as the argument for $L''$ is similar. The reader may find it helpful to refer to Figures \ref{fig:diag_sum_example} and \ref{fig:diagrammatic_sum_cartoon}. To show $L'$ is alternating, we must show that, tracing along $C$ in $L'$, for each $i\in [l]$, you exit the last crossing in $X_{(i+1)\mod l}$ via an understrand if and only if you encounter the first crossing in $X_i$ via an overstrand. This statement holds if and only if, tracing along $C$ in $L$, for each $i\in [l]$, you exit the last crossing in $X_{(i+1)\mod l}$ via an understrand if and only if you encounter the first crossing in $X_i$ via an overstrand. Fix $i\in[l]$. In $L$, the segment entering the first crossing of $X_i$ is the last segment exiting $Y_i$, and the segment exiting the last crossing of $X_{(i+1)\mod l}$ is the segment entering the first crossing of $Y_i$. Any pair of adjacent crossings in $Y_i$ share a segment, by construction. Since $L$ is alternating, the first segment entering $Y_i$ enters via an overstrand if and only if the last segment exiting $Y_i$ exits via an understrand. This implies the desired statement.  
    
    Finally, we address $\widetilde{L'}$ and again observe that the argument for $\widetilde{L''}$ is identical. The reader may find it helpful to refer to Figure \ref{fig:y_move}. Recall that because $L$ is alternating, so is $L'$. Additionally, we recall that $\widetilde{L'}$ is obtained from $L'$ via swap moves. As you trace along $C$ in $L'$, by the definition of Seifert circle and because $L'$ is alternating, you either always approach a group of crossings $X_i$ on the understrand and exit on the overstrand, or always approach on the overstrand and exit on the understrand. After performing a swap move to obtain $\widetilde{L'}$, the segment exiting the last crossing of $X_i$ becomes connected to the segment entering the first crossing of $X_i$. Thus we conclude that since $L'$ is alternating, so is $\widetilde{L'}$. 
 \end{proof}

\section{Proof of main result}
\label{sec:proof}

In the preceding section we studied the structure of graphs $\G_{L,i}$ for diagrams which are diagrammatic Murasugi sums $L=L'*_C L''$. In this section, we study their perfect matchings. Our goal is to partition the set of perfect matchings of the truncated face-crossing incidence graphs of alternating link diagrams in a way that will allow us to prove the main results of the paper. The key in establishing this partition is the following property.

\begin{lemma}\label{lem:even_flock}
Let $M$ be a perfect matching of $\G_{L,i}$. Then
\[\#\{e\in M_{\text{flock}} \ | \ e \text{ is incident to } r'_1, \dots, r'_{l-1}\} = \#\{e\in M_{\text{flock}} \ | \ e \text{ is incident to } r''_2, \dots, r''_l\},\]
where $M_{\text{flock}}$ is the set of edges of $M$ that are flock edges in $\G_{L,i}$. 
In particular, $M$ contains an even number of flock edges, namely, $2\#\{e\in M_{\text{flock}} \ | \ e \text{ is incident to } r'_1, \dots, r'_{l-1}\}.$
\end{lemma}

\begin{proof}
Let $M\in\mathcal{S}_{L,i}$. Since $\G_{L,i}$ is bipartite, for each $e\in M_{\text{flock}}$, $e$ meets precisely one black vertex and one white vertex of $\G_{L,i}$. Recall that, by definition, black vertices correspond to crossings of $L$ and white vertices correspond to regions. By Notation \ref{not:sums}, flock edges of $M$ satisfy the following property.
\begin{center}
    \textit{
    For each $e\in M_{\text{flock}}$, either $e$ is incident to $r'_j$ and a crossing in $Y_j$ for some $1\leq j\leq l-1$, or $e$ is incident to $r''_k$ and a crossing in $X_k$ for some $2\leq k \leq l$.
    }
\end{center}

Let $R'$ denote the set of white vertices of $L'$ incident to edges in $M_{\text{flock}}$, and let $C'$ denote the set of black vertices of $L'$ incident to edges in $M_{\text{flock}}$. By the above property,
\[\#\{e\in M_{\text{flock}} \ | \ e \text{ is incident to } r'_1, \dots, r'_{l-1}\} = \# R', \text{ and } 
 \#\{e\in M \ | \ e \text{ is incident to } r''_2, \dots, r''_l\} = \# C'.\]
 
 Suppose for contradiction that
\[\#\{e\in M_{\text{flock}} \ | \ e \text{ is incident to } r'_1, \dots, r'_{l-1}\} \neq \#\{e\in M \ | \ e \text{ is incident to } r''_2, \dots, r''_l\}.\]

Without loss of generality, 
\[\#\{e\in M_{\text{flock}} \ | \ e \text{ is incident to } r'_1, \dots, r'_{l-1}\} > \#\{e\in M \ | \ e \text{ is incident to } r''_2, \dots, r''_l\}.\]

By Lemma \ref{lem:d_move}, $M$ restricts to a matching of each of $\G_{L',i}$. Moreover, $M$ restricts to a perfect matching of $\G_{L',i}\setminus(R'\sqcup C')$. But, $\G_{L',i}\setminus(R'\sqcup C')$ is a bipartite graph with more black vertices than white vertices, and thus cannot admit a perfect matching. This is a contradiction.

\end{proof}

\begin{theorem}\label{thm:term_comparison}
    Let $L$ be a link diagram and suppose $L=L'*_C L''$ where $C$ has length $l \geq 2$. Let $\widetilde{L'}$ and $\widetilde{L''}$ be as in Definition \ref{def:y-moves-links}. Let $\mathcal{S}'_{L,i}$ denote the perfect matchings of $\G_{L,i}$ which contain $2k$ flock edges, where $1\leq k \leq l-2$,  of $\G_{L,i}$. Then
    \[\tDelta_L = \tDelta_{L'}\tDelta_{L''} + \tDelta_{\widetilde{L'}}\tDelta_{\widetilde{L''}} + \sum_{M \in \mathcal{S}'_{L,i}} \wt (M)\]

\end{theorem}

\begin{proof}
We show that we can group terms in the state-sum formula for $\tDelta_L$ to obtain the formula in the theorem statement. Recall that $\mathcal{S}_{L,i}$ denotes the set of perfect matchings of $\G_{L,i}$, and also recall that each $M\in\mathcal{S}_{L,i}$ contrubutes a term $\wt(M)$ to $\tDelta_L$. 

The matchings of $\G_{L,i}$ which use no flock edges are in weight-preserving bijection with matchings of $\G_{L',i'}\sqcup\G_{L'',i''}$. The matchings of $\G_{L,i}$ which use one flock edge incident to each $r'_j$, $j \leq l-1$, and each $r''_k$, $2\leq k \leq l$, are in weight-preserving bijection with matchings of $\G_{\widetilde{L'},\widetilde{i''}}\sqcup\G_{\widetilde{L'i},\widetilde{i''}}$. These bijections follow from Lemmas \ref{lem:d_move} and \ref{lem:y_moves}. 

Suppose $l>2$. By Lemma \ref{lem:even_flock}, an element of $\mathcal{S}_{L,i}$ must contain an even number of flock edges. Partitioning $\mathcal{S}_{L,i}$ into those perfect matching that (1) use $0$ flock edges (this is the minimum possible to use in a perfect matching of $\G_{L,i}$) and those that (2) use $2(l-1)$ flock edges (this is the maximum possible to use in a perfect matching of $\G_{L,i}$), and those that (3) use  $2k$ flock edges, where $1\leq k \leq l-2$, we obtain $\tDelta_L$ as the sum of $\tDelta_{L'}\tDelta_{L''}$ coming from (1), $\tDelta_{\widetilde{L'}}\tDelta_{\widetilde{L''}}$ coming from (2), and the weights of perfect matchings which use  $2k$ flock edges, where $1\leq k \leq l-2$ flock edges coming from (3).

Thus we obtain 
\begin{equation}
    \sum_{M \in \mathcal{S}_{L,i}} \wt(M) = \sum_{\substack{M' \in \mathcal{S}_{L',i'}\\M'' \in \mathcal{S}_{L'',i''}}} \wt(M') \wt(M'') + \sum_{\substack{N' \in \mathcal{S}_{\widetilde{L'},\widetilde{i'}}\\N'' \in \mathcal{S}_{\widetilde{L''},\widetilde{i''}}}} \wt(N') \wt(N'') +  \sum_{M \in \mathcal{S}'_{L,i}} \wt (M)
\end{equation}
which, using Theorem \ref{thm:dimer_version}, is the desired formula.

\end{proof}

When $C$ has length exactly 2, the last term in Theorem \ref{thm:term_comparison} is equal to zero, and we obtain the following result, which is a generalization of \cite[Proposition 2.21]{kalman}.

\begin{theorem}\label{thm:kalman_generalization} Suppose $L=L' *_C L''$ where $C$ has length 2. Let $\widetilde{L'}$ and $\widetilde{L''}$ be as in Definition \ref{def:y-moves-links}.
Then,

    \[\tDelta_L = \tDelta_{L'}\tDelta_{L''} + \tDelta_{\widetilde{L'}}\tDelta_{\widetilde{L''}}.\]
\end{theorem}

\begin{proof}

This follows from Lemma \ref{lem:even_flock} and Theorem \ref{thm:term_comparison}.

\end{proof}

We also collect the following results when $C$ has smaller length. 
\begin{remark}\label{rem:zero_length}
If $L= L' *_C L''$ where $C$ has length 0, then $L$ is disconnected and $\tDelta_L=0$. 
\end{remark}

\begin{lemma}\label{lem:length_1}
If $L= L' *_C L''$ where $C$ has length 1, then $\tDelta_L = \tDelta_{L'}\tDelta_{L''}$.
\end{lemma}

\begin{proof}
    The crossings incident to $C$ all meet an absent region of the diagram. Thus there are no flock edges in $\G_{L,i}$, and $\G_{L,i} = \G_{L', i'} \sqcup \G_{L'', i''}$ as weighted graphs. The state-sum formula implies that $\tDelta_L = \tDelta_{L'}\tDelta_{L''}$
\end{proof}

We now turn our attention to alternating links, where we will use the results above to prove certain $\tDelta_L$ are trapezoidal. The following notion will be useful to us.

\begin{definition}\label{def:non_cancel}
    Let $f=\sum_{i \in \frac{1}{2}\mathbb{Z}} a_i t^i$ and $g= \sum_{i \in \frac{1}{2}\mathbb{Z}} b_i t^i$ be Laurent polynomials in $t^{\pm 1/2}$. The sum $f+g$ is \textbf{non-canceling} if all terms of the same degree have the same sign; that is, for $i \in \frac{1}{2}\mathbb{Z}$, $a_i, b_i$ are either both nonnegative or both nonpositive. The product $fg$ is \textbf{non-canceling} if the sum $\sum_{i} ( \sum_{j+k=i} a_j b_k t^i)= fg$ is non-canceling; that is, for each $i$, either $a_j b_k$ is nonnegative for all $j,k$ with $i=j+k$, or $a_j b_k$ is nonpositive for all $j,k$ with $i=j+k$.
\end{definition}

We will need the following lemma, which shows that the state-sum formula for $\tDelta_L$ is non-canceling if $L$ is alternating. It follows readily from Kauffman's Clock Theorem \cite[Theorem 2.5]{K06}.

\begin{lemma}\label{rem:non_cancel}
Suppose $L$ is alternating. Then if $M, M' \in \mathcal{S}_{L,i}$ are such that $\wt(M)$ and $\wt(M')$ are the same degree, then $\wt(M)$ and $\wt(M')$ are also the of same sign. In particular, the sum in Theorem \ref{thm:dimer_version} is non-canceling.
\end{lemma}

\begin{proof}
    Kauffman's Clock Theorem shows that any two states of a connected link diagram are related by a sequence of \emph{clock moves}. Each clock move changes the sign of the weight. If the link diagram is alternating, then each clock move also changes the degree by 1. The lemma follows.
\end{proof}

We next will show that for alternating diagrams $L = L' *_C L''$, the polynomials $\tDelta_{L'}\tDelta_{L''}$ and $\tDelta_{\widetilde{L'}}\tDelta_{\widetilde{L''}}$ ``fit inside" the Alexander polynomial $\tDelta_L$ with no cancellation. This requires a bit of setup.

\begin{notation}\label{not:abs_val} For 
\[f=\sum_{i \in \frac{1}{2} \mathbb{Z}} a_i t^i \in \mathbb{Z}[t^{\pm 1/2}]\]
we use the notation $|f|:= \sum_i|a_i| t^i$.
\end{notation}

\begin{lemma}\label{lem:abs_val_sum_prod}
Let $f=\sum_{i \in \frac{1}{2}\mathbb{Z}} a_i t^i$ and $g= \sum_{i \in \frac{1}{2}\mathbb{Z}} b_i t^i$ be Laurent polynomials in $t^{\pm 1/2}$. If the sum $f+g$ is non-canceling, then $|f+g| = |f| + |g|$. If the product $fg$ is non-canceling, then $|fg|=|f||g|$.
\end{lemma}

\begin{proof}
    The statement about the sum is immediate from the definitions.

    For the statement about the product, by assumption $fg = \sum_{i} ( \sum_{j+k=i} a_j b_k t^i) $ is non-canceling. This means $|fg| =\sum_{i} ( \sum_{j+k=i} |a_j b_k| t^i) $, which is clearly equal to $|f||g|$.
\end{proof}

\begin{theorem}\label{thm:alternating-sums-non-canceling}
    Suppose $L$ is alternating and $L= L' *_C L''$. If $C$ has length $1$, then $|\tDelta_L| =  |\tDelta_{L'}||\tDelta_{L''}|$. If $C$ has length $l \geq 2$, then the right-hand sides of Theorem \ref{thm:term_comparison}, Theorem \ref{thm:kalman_generalization}, and Lemma \ref{lem:length_1} are non-canceling. 
    
    In particular, for an alternating link $L$ we have that, coefficient-wise,
    $$|\tDelta_L| \geq  |\tDelta_{L'}\tDelta_{L''} + \tDelta_{\widetilde{L'}}\tDelta_{\widetilde{L''}}| = |\tDelta_{L'}||\tDelta_{L''}| + |\tDelta_{\widetilde{L'}}||\tDelta_{\widetilde{L''}}| .$$
   
\end{theorem}

\begin{proof} If $C$ has length $1$, the result follows by Lemma \ref{lem:length_1}. Suppose $C$ has length $l \geq 2$. The equality in Theorem \ref{thm:term_comparison} is proven via a weight-preserving bijection of matchings. This shows that the right-hand side of the equality 
\[\tDelta_L = \tDelta_{L'}\tDelta_{L''} + \tDelta_{\widetilde{L'}}\tDelta_{\widetilde{L''}} + \sum_{M \in \mathcal{S}'_{L,i}} \wt (M)\] is non-canceling. Thus, by Lemma \ref{lem:abs_val_sum_prod},

\begin{eqnarray*}|\tDelta_L| &= |\tDelta_{L'}\tDelta_{L''} + \tDelta_{\widetilde{L'}}\tDelta_{\widetilde{L''}} + \sum_{M \in \mathcal{S}'_{L,i}} \wt (M)| \\ 
 &= |\tDelta_{L'}\tDelta_{L''} + \tDelta_{\widetilde{L'}}\tDelta_{\widetilde{L''}}| + |\sum_{M \in \mathcal{S}'_{L,i}} \wt (M)| \\
 &\geq |\tDelta_{L'}\tDelta_{L''} + \tDelta_{\widetilde{L'}}\tDelta_{\widetilde{L''}}|, 
\end{eqnarray*}
where the last inequality is considered coefficientwise.

Also, by Lemma \ref{lem:abs_val_sum_prod},
\[|\tDelta_{L'}\tDelta_{L''} + \tDelta_{\widetilde{L'}}\tDelta_{\widetilde{L''}}| = |\tDelta_{L'}||\tDelta_{L''}| + |\tDelta_{\widetilde{L'}}||\tDelta_{\widetilde{L''}}| .\]
\end{proof}

The next result shows that $\tDelta_{L'}\tDelta_{L''}$ covers the entire ``spread" of $\tDelta_L$.    Recall that   the \textit{support} of  $f\in\mathbb{Z}[t^{1/2},t^{-1/2}]$ is the subset $S \subset \frac{1}{2}\mathbb{Z}$ such that the monomial $t^s$ appears in $f$ with nonzero coefficient if and only if $s \in S$. We define the \emph{degree} of $f\in\mathbb{Z}[t^{1/2},t^{-1/2}]$ to be (max element of the support)$-$(min element of the support).

\begin{theorem} \label{thm:degree}
     Suppose $L$ is an alternating link diagram and $L = L' *_C L''$ where $C$ has length $l$. Then $\deg(\tDelta_{L'}\tDelta_{L''})=\deg(\tDelta_L)$ and, as long as $\tDelta_{\widetilde{L'}}\tDelta_{\widetilde{L''}}$ is nonzero, $\deg(\tDelta_{\widetilde{L'}}\tDelta_{\widetilde{L''}})=\deg(\tDelta_L)-2l+2$.
\end{theorem}

\begin{proof}
    Let $d$, $d'$, $d''$, $\widetilde{d'}$, and $\widetilde{d''}$ denote the number of crossings of $L$, $L'$, $L''$, $\widetilde{L'}$, and $\widetilde{L''}$, respectively. Let $f$, $f'$, $f''$,  $\widetilde{f'}$, and $\widetilde{f''}$ denote the number of Seifert cycles of $L$, $L'$, $L''$, $\widetilde{L'}$, and $\widetilde{L''}$, respectively. In \cite{crowell}, the author proves that for a connected\footnote{A link diagram $L$ is connected if there does not exist a closed path in the plane separating $L$ into two nontrivial diagrams.} diagram of an alternating link, the degree of the Alexander polynomial is $d - f +1$. If $L$ is a disconnected diagram, then its Alexander polynomial is zero.

    By construction of $L'$ and $L''$, we have $d=d'+d''$. Also in this construction, the Seifert circle $C$ of $L$ contributes one Seifert circle to each of $L'$ and $L''$. Call these newly created Seifert cycles $C'$ and $C''$, respectively. All other Seifert circles of $L$ become a Seifert circle of either $L'$ or $L''$, so $f= f' + f'' -1$. If $L$ is disconnected, then so are so are $L'$ and $L''$, meaning $\tDelta_L$ and $\tDelta_{L'}\tDelta_{L''}$ are both zero. If $L$ is connected, then so are $L'$ and $L''$, and the degree of $\tDelta_{L'}\tDelta_{L''}$ is 
    $$(d' - f' +1)+ (d''-f'' +1) =d' + d'' - (f' + f'' -1) +1 = d -f +1.  $$
    So $\tDelta_L$ and $\tDelta_{L'}\tDelta_{L''}$ have the same degree. 
    
    Recall, in the construction of $\widetilde{L'}$ from $L'$, we have $d'=\widetilde{d'}$. All Seifert circles of $L'$ which are not $C'$ are present in $\widetilde{L'}$ as well. However, the swap move replaces $C'$ with $l$ distinct Seifert circles. That is, $\widetilde{L'}$ has $l-1$ more Seifert circles than $L'$.  If either $\widetilde{L'}$ or $\widetilde{L''}$ is disconnected, then $\tDelta_{\widetilde{L'}}\tDelta_{\widetilde{L''}}=0$. Suppose both $\widetilde{L'}$ and $\widetilde{L''}$ are connected link diagrams. Then, the degree of $\tDelta_{\widetilde{L'}}$ is $\widetilde{d'}-\widetilde{f'}+1=d'-(f'+(l-1))+1$. Similarly, after applying the swap moves, $C''$ contributes $l$ Seifert circles to $\widetilde{L''}$. So, the degree of $\tDelta_{\widetilde{L''}}$ is $d''-(f''+(l-1))+1$. Therefore, the degree of $\tDelta_{\widetilde{L'}}\tDelta_{\widetilde{L''}}$ is 
    $$(d' - (f' + (l-1)) +1)+ (d''-(f''+(l-1)) +1) = (d' + d'' - (f' + f'' -1) + 1) - 2l  + 2 = (d-f+1) - 2l  + 2.$$
    
 This proves the result.
    
\end{proof}

From this, we readily prove the following result.

\begin{theorem} \label{thm:terms}
     Suppose $L$ is an alternating link diagram and $L = L' *_C L''$. Then $\tDelta_L$ and $\tDelta_{L'}\tDelta_{L''}$ have the same support.    
\end{theorem}

\begin{proof}
For an alternating link $L$, the number of terms in its Alexander polynomial $\tDelta_L$ is one more than the degree, as a consequence of the Clock Theorem. Since the sum
    \[\tDelta_{L'}\tDelta_{L''} = \sum_{M \in \mathcal{S}_{L,i}} \wt(M) = \sum_{\substack{M' \in \mathcal{S}_{L',i'}\\M'' \in \mathcal{S}_{L'',i''}}} \wt(M') \wt(M'')\]
    is non-canceling, the number of terms in $\tDelta_{L'}\tDelta_{L''}$ is also one more than its degree. By Theorem \ref{thm:degree}, the degrees of $\tDelta_L$ and $\tDelta_{L'}\tDelta_{L''}$ are the same, and therefore, the number of terms in both polynomials is the same. Moreover, since the supports of both $\tDelta_L$ and $\tDelta_{L'}\tDelta_{L''}$ are centered around zero, it follows  that their supports are the same.
    
\end{proof}

Now, we turn to trapezoidal properties. We need the following lemmas.

\begin{lemma}\label{lem:murasugi_trapezoidal}\cite[Proposition 2.1]{murasugi}
Let $A(x) = \sum_{i=0}^n a_ix^i$, $B(x) = \sum_{i=0}^m b_ix^i$ be polynomials with positive coefficients such that $a_0, \dots, a_n$ and $b_0, \dots, b_m$ are both trapezoidal sequences. Then the coefficient sequence of $A(x)B(x)$ is also trapezoidal.
\end{lemma}

\begin{definition}
    We say that a Laurent polynomial $f=\sum_i a_it^i \in\mathbb{Z}[t^{1/2},t^{-1/2}]$ is \textbf{centered around zero} if $f= \pm \sum_i a_it^{-i}.$
\end{definition}

\begin{lemma}\label{lem:symm_sum_and_prod}

If $A(t), B(t) \in \mathbb{Z}[t^{1/2},t^{-1/2}]$ are centered around zero, then so are $A(t)+B(t)$ and $A(t)B(t)$.

\end{lemma}

\begin{lemma}\label{lem:trap_sum_and_prod}
If $f,g\in\mathbb{Z}[t^{1/2},t^{-1/2}]$ are trapezoidal Laurent polynomials centered around zero, then $|f||g|$ is trapezoidal and centered around zero. If additionally $\deg(f) = \deg(g) + 2$, then $|f| + |g|$ is trapezoidal and centered around zero.
\end{lemma}
\begin{proof}
    The statement about the product is immediate from Lemmas \ref{lem:murasugi_trapezoidal} and \ref{lem:symm_sum_and_prod}. For the statement about the sum, being centered around zero is immediate from \cref{lem:symm_sum_and_prod}. The degree condition implies that either all exponents of $f,g$ are integral, or all exponents of $f,g$ are half-integral, so the same is true of $|f|+|g|$. Clearly the coefficients of $|f|+|g|$ are unimodal (that is, weakly increasing, and then weakly decreasing). That they are also trapezoidal, meaning that they have at most one ``plateau" follows from the degree condition.
\end{proof}

We prove in Theorem \ref{thm:full_decomp} below that as long as the \ttwo circles of $L$ are small in length, $|\tDelta_L|$ is trapezoidal. This result is originally due to Azarpendar, Juh\'asz and K\'alm\'an    \cite[Theorem 2.23]{kalman}. Azarpendar, Juh\'asz and K\'alm\'an  \cite{kalman}   proved  the case where   each \ttwo cycle $C$, $C$ has length two and $|X_1|,|X_2|,|Y_1|,|Y_2|=1$. They present a proof sketch for the general case of the theorem. Below, we prove the theorem in full generality with all details included. The viewpoint of perfect matchings streamlines the proof considerably.  Our inductive method of the following theorem originates from the proof of \cite[Theorem 2.22]{kalman}, where the authors use an induction for the aforementioned specialized statement. We induct on the number of \ttwo cycles of $L$; the authors of \cite{kalman} induct on the number of terms in the diagrammatic Murasugi sum of $L$ into special alternating links, which is the number of \ttwo cycles plus one.  

\begin{theorem}\label{thm:full_decomp} \cite[Theorem 2.23]{kalman}
    Suppose $L$ is an alternating link diagram such that all \ttwo circles of $L$ are length at most 2. Then $|\tDelta_L|$ is trapezoidal.
\end{theorem}

\begin{proof}
We induct on the number of \ttwo circles of $L$. If $L$ has no \ttwo circles, it is a special alternating diagram, and $|\tDelta_L|$ is trapezoidal by \cite[Theorem 1.2]{HMV-paper1}. 

Now, suppose $C$ is a \ttwo circle of $L$. If $C$ has length 0, then $L$ is disconnected, so $|\tDelta_L|=0$ and in particular is trapezoidal. If $C$ has positive length, we have $L= L' *_C L''$.

\noindent \textbf{Case 1:} If $C$ has length 1, then by Lemma \ref{lem:length_1}, $\tDelta_L = \tDelta_{L'}\tDelta_{L''}$. By Remark \ref{rem:non_cancel}, right-hand side is non-canceling, so by Lemma \ref{lem:abs_val_sum_prod}, $$|\tDelta_L| = |\tDelta_{L'}|~|\tDelta_{L''}|.$$  

By Remark \ref{rem:decomp_tone} and Lemma \ref{lem:y-moves-alternating-and-type-2-circles}, we have that $L', L''$ are alternating and have fewer \ttwo circles than $L$, which are again all of length at most 2 as they are inherited from $L$. Thus by induction hypothesis $|\tDelta_{L'}|$ and $|\tDelta_{L''}|$ are trapezoidal. By Lemma \ref{lem:trap_sum_and_prod}, their product $|\tDelta_{L'}|\cdot |\tDelta_{L''}| = |\tDelta_L|$ is also trapezoidal, completing the inductive step. 

 \noindent \textbf{Case 2:} If $C$ has length $2$, then by Theorem \ref{thm:kalman_generalization}, 
$$\tDelta_L = \tDelta_{L'}\tDelta_{L''} + \tDelta_{\widetilde{L'}}\tDelta_{\widetilde{L''}}$$
 and by Theorem \ref{thm:alternating-sums-non-canceling}, the right-hand side is non-canceling. So by Lemma \ref{lem:abs_val_sum_prod}, we also have 
 $$|\tDelta_L| = |\tDelta_{L'}|~|\tDelta_{L''}| + |\tDelta_{\widetilde{L'}}|~|\tDelta_{\widetilde{L''}}|.$$

By Remark \ref{rem:decomp_tone} and Lemma \ref{lem:y-moves-alternating-and-type-2-circles}, $L'$, $L''$, $\widetilde{L'}$, and $\widetilde{L''}$ are alternating and have fewer \ttwo circles than $L$. Again all of these \ttwo circles are length at most 2, because they are inherited from $L$. By induction hypothesis $|\tDelta_{L'}|$, $|\tDelta_{L''}|$, $|\tDelta_{\widetilde{L'}}|$, and $|\tDelta_{\widetilde{L''}}|$ are trapezoidal. Thus, by Lemma \ref{lem:trap_sum_and_prod}, the two products $|\tDelta_{L'}|\cdot |\tDelta_{L''}|$ and $|\tDelta_{\widetilde{L'}}|\cdot |\tDelta_{\widetilde{L''}}|$ are also trapezoidal and centered around zero. By \cref{thm:degree}, the degrees of the two terms differ by $-2l+2 = -2$, provided that both terms are nonzero. So using Lemma \ref{lem:trap_sum_and_prod} again, their sum is trapezoidal.
\end{proof}

\section{Obstacles to extending {Theorem \ref{thm:full_decomp}} past length 2}

\label{sec:obstacle}
\subsection{Obstacles for using Theorem \ref{thm:term_comparison}}
Consider the formula from Theorem \ref{thm:term_comparison}:
    \begin{equation}\label{eq:break-down-alexander-poly}
    \tDelta_L = \tDelta_{L'}\tDelta_{L''} + \tDelta_{\widetilde{L'}}\tDelta_{\widetilde{L''}} + \sum_{M \in \mathcal{S}'_{L,i}} \wt (M)\end{equation}
where $\mathcal{S}'_{L,i}$ denotes perfect  matchings which use $2k$ flock edges for $1\leq k \leq l-2$.

Equation \eqref{eq:break-down-alexander-poly} suggests an inductive proof strategy for showing that $|\tDelta_L|$ is trapezoidal. The ``controllable" terms 
\[\tDelta_{L'}\tDelta_{L''} \qquad \text{and} \qquad \tDelta_{\widetilde{L'}}\tDelta_{\widetilde{L''}}\]
are each sums over matchings of smaller (truncated face-crossing incidence) graphs. If the ``remainder term"
\[\sum_{M \in \mathcal{S}'_{L,i}} \wt (M)\]
could also be organized into sums over matchings of smaller graphs, one could hope to understand those individual sums inductively as well. For the induction to work, each individual sum should be a trapezoidal polynomial centered around zero.

\begin{question}\label{goal:partition-remainder-term} Can $\mathcal{S}'_{L,i}$ be partitioned so that each block consists of matchings of a subgraph of $\G_{L,i}$, and summing over each block gives a centered trapezoidal polynomial? 
\end{question}

\begin{figure}
    \centering
    \includegraphics[width=0.9\textwidth]{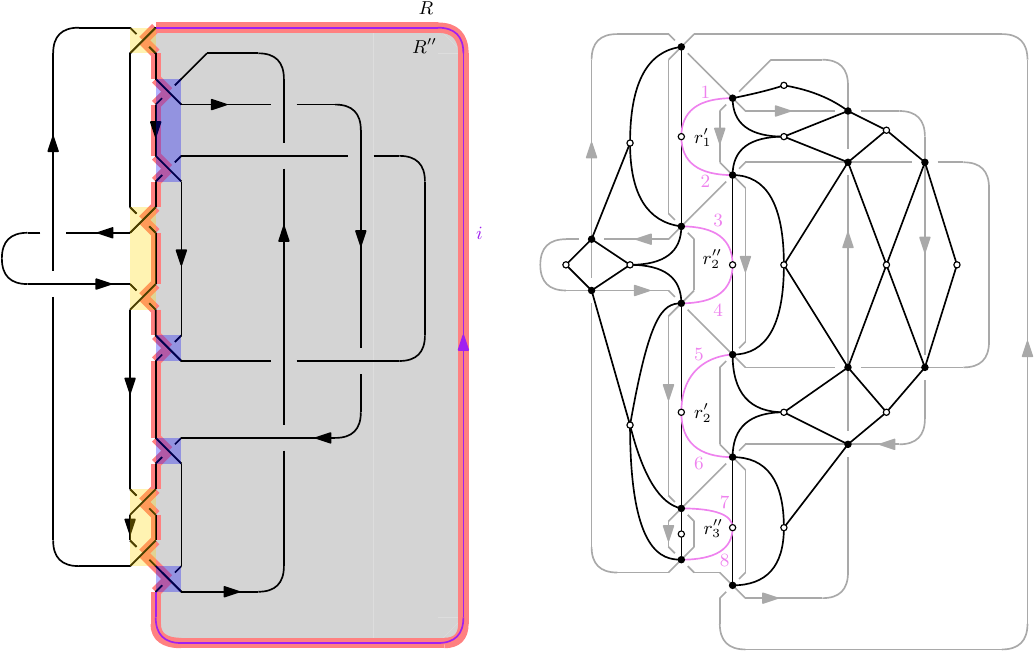}
    \caption{An alternating diagram $L$ which is the length-3 sum of two special alternating links and its truncated face-crossing incidence graph. On the left, the \ttwo circle is in red, the yellow crossings belong to $R'$ and the blue crossings belong to $R''$. On the right, the flock edges are in pink.}
    \label{fig:length3}
\end{figure}

The most natural way to partition $\mathcal{S}'_{L,i}$ is according to which flock-edges are used in the matching. We analyze this partition for the example in Figure \ref{fig:length3}, which has 
\[|\tDelta_L|= 14 t^{-9/2}+108t^{-7/2}+395t^{-5/2}+882t^{-3/2} + 1320t^{-1/2}+1302 t^{1/2}+882 t^{3/2}+395 t^{5/2}+108 t^{7/2}+14 t^{9/2}.\] Unfortunately, as we will see, partitioning $\mathcal{S}'_{L,i}$ according to flock edges used is not compatible with Question \ref{goal:partition-remainder-term}.

Figure \ref{fig:length3} shows an alternating length 3 Murasugi sum of two special alternating diagrams. Matchings of $\G_{L,i}$ use either 0, 2, or 4 flock edges, by Lemma \ref{lem:even_flock}. So $\mathcal{S}'_{L,i}$ consists of the matchings which use exactly 2 flock edges. For such matchings, one flock edge must be incident to $r_i'$ for some $i$ and the other must be incident to $r_j''$ for some $j$, again by Lemma \ref{lem:even_flock}. So there are 16 possible pairs of flock edges that can be used in $M \in \mathcal{S}'_{L,i}$, listed in the left column of Table \ref{tab:two-flock-edge-contributions}. The right column gives the sum of $|\wt(M)|$ for matchings using exactly that pair of flock edges.

\begin{table}[htbp!]
\renewcommand{\arraystretch}{1.3}
    \centering
    \begin{tabular}{|c|c|cccccccc|}
    \hline
     &Edges in $M$ & \multicolumn{8}{c|}{Contribution to $|\tDelta_L|$}  \\ 
     \hline
      A  & 1,3    &$4t^{-7/2} +$&$ 17t^{-5/2} +$&$ 35t^{-3/2} +$&$ 49t^{-1/2} + $&$47 t^{1/2} +$&$ 30 t^{3/2} +$&$ 12 t^{5/2} + $&$2 t^{7/2} $\\ \hline
 
  B&  1,4     & & & $4t^{-3/2}+$&$ 9t^{-1/2}+$&$ 9 t^{1/2} +$&$ 5 t^{3/2} +$&$ t^{5/2}$ &\\ \hline
  
C     & 1,7   && $2t^{-5/2} +$&$ 2t^{-3/2} +$&$t^{-1/2}+ $&$t^{1/2}$&&& \\ \hline
  
  D    & 1,8   &  $2t^{-7/2} +$&$ 2t^{-5/2} +$&$ t^{-3/2} +$&$t^{-1/2}$&&&&\\ \hline

E    & 2,3   &$4t^{-7/2} +$&$ 23t^{-5/2} +$&$ 62t^{-3/2} +$&$ 105t^{-1/2} +$&$ 123 t^{1/2} +$&$ 100 t^{3/2} +$&$ 52 t^{5/2} +$&$ 14 t^{7/2}$ \\ \hline

F  &  2,4  && &$4t^{-3/2} +$&$ 15t^{-1/2} +$&$ 24 t^{1/2} +$&$ 19 t^{3/2} +$&$ 7 t^{5/2} $ &\\ \hline

G     & 2,7  & & $2t^{-5/2} +$&$ 5t^{-3/2} +$&$3t^{-1/2} $ & & & &\\ \hline

H     & 2,8   & $2t^{-7/2} +$&$ 5t^{-5/2} +$&$ 3t^{-3/2}$ & & & & & \\ \hline

I     & 3,5  & & $7t^{-5/2} +$&$ 26t^{-3/2}+$&$ 41t^{-1/2} +$&$ 35 t^{1/2} +$&$ 17 t^{3/2} +$&$ 4 t^{5/2}$ & \\ \hline

J     & 3,6  & & $t^{-5/2} +$&$ 5t^{-3/2}+$&$11 t^{-1/2}+ $&$11 t^{1/2} +$&$ 4 t^{3/2}$& &\\ \hline

K     & 4,5   &  $14t^{-7/2} +$&$ 66t^{-5/2} + $&$141t^{-3/2} +$&$ 178t^{-1/2}+ $&$145 t^{1/2} +$&$ 77 t^{3/2} + $&$25 t^{5/2} + $&$4 t^{7/2}$ \\ \hline

L     & 4,6  & $2t^{-7/2} +$&$ 12t^{-5/2} +$&$ 33t^{-3/2} +$&$49t^{-1/2}+$&$ 41 t^{1/2} +$&$ 19 t^{3/2} +$&$ 4 t^{5/2}$ & \\ \hline

M     & 5,7 & & & $4t^{-3/2}+$&$16t^{-1/2} + $&$25 t^{1/2} +$&$ 19 t^{3/2} +$&$ 6 t^{5/2}$ &\\ \hline

N     & 5,8 & &$4t^{-5/2} +$&$ 16t^{-3/2} +$&$25t^{-1/2}+ $&$19 t^{1/2} + $&$6 t^{3/2}$ & &\\  \hline

O     & 6,7  & & $6t^{-5/2} +$&$ 31t^{-3/2} +$&$73t^{-1/2}+ $&$105 t^{1/2} +$&$ 96 t^{3/2} + $&$53 t^{5/2} +$&$ 14 t^{7/2}$ \\ \hline

P     & 6,8  & $6t^{-7/2} +$&$ 31t^{-5/2} +$&$ 73t^{-3/2}+$&$105t^{-1/2} +$&$ 96 t^{1/2} +$&$ 53 t^{3/2} + $&$14 t^{5/2}$ &\\  \hline
    \end{tabular}
    \caption{Each line gives the contribution to $|\tDelta_L|$ from the matchings containing exactly the two flock edges indicated.}
    \label{tab:two-flock-edge-contributions}
\end{table}

The following sums of polynomials in Table \ref{tab:two-flock-edge-contributions} are centered and trapezoidal:
\[B + J;\quad G + M + N;\quad G + I + J + M; \quad B + C + G + I + M;\quad B + G + J + M + N; \quad B + E +F +G + K+ M; \]
\[ B + H + K + L + M + O; \quad A + B + H + K + I + M + O\]
as well as the complementary sums. However, none of these sums are sums over matchings of subgraphs of $\G_{L,i}$. So, partitioning $\mathcal{S}'_{L,i}$ according to which flock edges are used in each matching does not help solve Question \ref{goal:partition-remainder-term}.

It is important to note, that in the handful of  examples we looked at, the entire polynomial $|\tDelta_L|-|\tDelta_{L'}||\tDelta_{L''}|- |\tDelta_{\widetilde{L'}}||\tDelta_{\widetilde{L''}}|$ was centered and trapezoidal for alternating links $L=L' *_C L''$ (and in particular, is centered and trapezoidal for this example). If one could prove this in general, it would settle Fox's conjecture.

\subsection{Obstacles for extending Theorem \ref{thm:kalman_generalization}}
 Theorem \ref{thm:kalman_generalization} deals with length 2 Murasugi sums, and one may hope to extend the result to sums of longer length. The proof can be rephrased as follows: partition the set of perfect matchings $\mathcal{S}_{L,i}$ into those which use the top flock edge of $\G_{L,i}$ and those which do not. Call these partition elements $\mathcal{N}$ and $\mathcal{N'}$, respectively. Then, analyze the sums over the weights of matchings in each of these two subsets of $\mathcal{S}_{L,i}$. These turn out to be centered trapezoidal polynomials, thereby proving the theorem.

For a Murasugi sum of longer length, the polynomials $\sum_{M\in\mathcal{N}}\wt(M)$ and $\sum_{M\in\mathcal{N'}}\wt(M)$ may not be centered nor trapezoidal. But one can generalize this partition by forming a binary tree $\mathcal{T}$ of subsets of $\mathcal{S}_{L,i}$. The root is $\mathcal{S}_{L,i}$. The two children of a node $\mathcal{M}$ at level $j$ are the subsets $\mathcal{N} = \{M \in \mathcal{M} : M ~\text{uses flock edge $j$}\}$ and $\mathcal{N'} = \{M \in \mathcal{M} : M ~\text{does not use flock edge $j$}\}$. The nodes at a fixed level of $\mathcal{T}$ give a partition of $\mathcal{S}_{L,i}$. The following question seeks to generalize the proof of Theorem \ref{thm:kalman_generalization}.

\begin{question}\label{quest:tree}
    For $L$ an alternating link with a length $>2$ \ttwo circle, is there a level of $\mathcal{T}$ such that, for each node $\mathcal{M}$ at that level, $\sum_{M\in\mathcal{M}}\wt(M)$ is centered and trapezoidal?
\end{question}

The example in Figure \ref{fig:length3} shows that the answer to Question \ref{quest:tree} is ``no." Only 8 non-root nodes of the tree give centered trapezoidal polynomials. These nodes are
\[\text{level 6: }000000, 101001;\quad \text{ level 7: } 0000000,1010010,1010011; \quad \text{level 8: }00000000,10100101,10100110 \]
where we use e.g. $101001$ to denote the subset of matchings which use flock edges 1, 3, 6, do not use flock edges 2, 4, 5 and may or may not use flock edges 7,8. 
No subset of these nodes comprises a partition of $\mathcal{S}_{L,i}$, so $\mathcal{T}$ is not useful in decomposing $|\tDelta_L|$ as a sum of centered trapezoidal polynomials.

\section*{Acknowledgements}   The authors thank David Cimasoni for his careful reading and interest in our work. KM is supported by the National Science Foundation under Award No.~DMS-2348676.
 MSB is supported by the National Science Foundation under Award No.~DMS-2444020. 
  Any opinions, findings, and conclusions or recommendations expressed in this material are
those of the author(s) and do not necessarily reflect the views of the National Science
Foundation.

\bibliographystyle{alpha}
\bibliography{refs}

\end{document}